\documentstyle[12pt,leqno]{article}
\begin{document}
\baselineskip 0.585cm
\def \bull{\vrule height .9ex width .8ex depth -.1ex}
\def \R{I \!\! R}
\def \N{I \!\! N}
\def \m{\mu}
\def \mac{\mu _{a,c}}
\def \g{\gamma}
\def \G{\Gamma}
\def \a{\alpha}
\def \k{\kappa}
\def \kac{\kappa _{a,c}}
\def \gs{\gamma^*}
\def \ah{\hat a}
\def \as{\alpha^*}
\def \at{\tilde \alpha}
\def \hg{H^{(\gamma) }}
\def \sh{\hat S}
\def \kg{K^{(\gamma)}}
\def \cag{C^{(\alpha,\gamma)}}
\def \bag{B^{(\alpha,\gamma)}}
\def \cags{C_l^{(\alpha^*_l,\gamma^*_l)}}
\def \bags{B_l^{(\alpha^*_l,\gamma^*_l)}}
\def \xg{\xi_\gamma}
\def \xga{\xi_{\alpha,\gamma}}
\def \xgn{\xi_{\gamma,n}}
\def \xgam{\xi_{\a,\gamma,2m}}
\def \gsu{generalized sieved ultraspherical~}
\def \nn{\nonumber}
\def \ms{\mu ^*}
\parindent 0.0cm
\begin{center}
  {\Large \bf Characterizations
  of generalized Hermite and sieved ultraspherical
  polynomials}\\

\bigskip
\bigskip

  {\sc By }  HOLGER DETTE\\
  {\it Institut f\"ur Mathematische Stochastik, Technische 
  Universit\"{a}t Dresden,\\
Mommsenstr. 13, 01062 Dresden, Germany.}\\

\bigskip
\medskip

  {\sc Abstract}\\
\end{center}

\medskip

A new characterization  of the generalized Hermite polynomials and
of the orthogonal polynomials with respect to the maesure $|x|^\g
(1-x^2)^{\a-1/2}dx$ is derived 
 which is based on a "reversing property" of the
coefficients in the corresponding recurrence formulas and does
not use the representation in terms of generalized Laguerre and Jacobi
polynomials. A similar characterization
can be obtained for a generalization of the sieved ultraspherical 
polynomials of the first and second kind. These results are applied
in order to determine  the asymptotic limit distribution for the zeros
when the degree and the parameters tend to infinity with the same
order.

\medskip

\smallskip


\bigskip

\begin{center}
{\sc 1. Introduction}\\
\end{center}
\def\theequation{1^.\arabic{equation}}
\setcounter{equation}{0}

Consider the generalized Hermite polynomials $\hg _n(x)$
orthogonal with respect to the measure $|x|^\g\exp (-x^2)dx$
($\g >-1$).
A characterization of these polynomials can easily be obtained from a
characterization of the generalized Laguerre polynomials
(see e.g. Al-Salam [3] or Chihara  [6]) and the well known 
relations between Laguerre- and Hermite polynomials (see [6,
p. 156]). In this paper we present two new characterizing properties
of the generalized Hermite polynomials which are based on a "reversing property"
of the coefficients in the corresponding three term recurrence 
relations and do not use explicitly the relation to the generalized
Laguerre polynomials. Similar
results can be derived for orthogonal polynomials
on a compact interval, say $[-1,1]$. Here the analogue of
the generalized Hermite polynomials (with respect to the new
characterizations) are the orthogonal polynomials with respect
to the measure $|x|^\g (1-x^2)^{\a-1/2}dx$ (see [6, p. 156])
which  satisfy a similar "reversing property". While these
"semi" classical orthogonal polynomials are discussed in Section 2,
Section 3 investigates analogous
 characterizations  for some "relatively new" systems of orthogonal
polynomials, namely the sieved random walk polynomials 
introduced by Charris and Ismail [4] (see also
Ismail [10], Charris and Ismail [5] and Geronimo and VanAssche
[9]). We give a characterization for a generalization of 
the sieved ultraspherical polynomials of the first and second 
kind,  $\cag _n(x,k)$, $ \bag _n(x,k)$, which are orthogonal
with respect to the measures 
$(1-x^2)^{\a-1/2}|U_{k-1}(x)|^{2\a}|T_k(x)|^\g$ and
$(1-x^2)^{\a+1/2}|U_{k-1}(x)|^{2\a}|T_k(x)|^\g$
(here $k \in \N$ is a fixed integer 
and $T_k(x)$ and $U_{k-1}(x)$
denote the Chebyshev polynomials of the first and second kind).
In the case $\g =0$ these polynomials were introduced by 
Al-Salam, Allaway and Askey [2]
as a limit from the $q-$ultraspherical polynomials of Rogers [12].

Finally the results of Section 2 and 3 are applied in Section 4 in order to derive the
asymptotic zero distribution of the polynomials $\hg _n(x)$,
$\cag _n(x,k)$ and $\bag _n(x,k)$ when the degree and parameters tend
 to infinity with the same order. For the generalized Hermite polynomials
 we thus  obtain an alternative proof of a recent result in
 Gawronski [8]
 while the limit distribution for the zeros of the \gsu 
 polynomials of the first and second kind is a generalization of 
 the measure for which 
 the sieved Chebyshev polynomials of the first
 kind are orthogonal with respect to (see [10])
 
 \bigskip

\begin{center}
  {\sc 2. Generalized Hermite polynomials and the weight function 
  $|x|^{\g}(1-x^2)^{\a-1/2}$}\\
\end{center}
\def\theequation{2.\arabic{equation}}
\setcounter{equation}{0}
Let $\m$ denote a probability measure on the real line $(-\infty,\infty)$
with all moments existing. The Stieltjes transform of $\m$
has the continued fraction expansion
\begin{equation}
\Phi (z)~= ~\int_{-\infty}^\infty {d\mu (x)\over z - x} ~=~
{~~~1~~|\over |z-b_1}~-~{~~a_1~~|\over|z-b_2}~-
~{~~a_2~~|\over|z-b_3}~-~\ldots
\label{cfg}
\end{equation}
where the quantities $a_i \geq 0$, $b_i \in \R$ ($i\geq 1$)
can be expressed in terms of the ordinary moments of $\m$
(see Perron [11] or Wall [15]). In the following
we consider the $n$th "terminated" continued fraction
of (\ref{cfg}) (i.e. we put $a_{n+1}=0$) and its corresponding
probability measure $\m_n$, that is 
\begin{equation}
\Phi_n (z)~= ~\int_{-\infty}^\infty {d\mu _n (x)\over z - x} ~=~
{~~1~~~|\over |z-b_1}~-~{~~a_1~~|\over|z-b_2}~-~ \ldots ~-~
{~~~~a_n~~~|\over|z-b_{n+1}} 
 \label{cft}
\end{equation}
It is well known that $\m_n$ has finite support given by the zeros of the
polynomial in the denominator of (\ref{cft}). Moreover, it can be shown (see
Dette and Studden [7, p. 4]) that the measure
$\m^R_n$ with Stieltjes transform corresponding to
the "reversed" continued fraction
\begin{equation}
\Phi^R_n(z)~= ~\int_{-\infty}^\infty {d\mu^R_n (x)\over z - x} ~=~
{~~~~1~~~~|\over |z-b_{n+1}}~-~{~~a_n~~|\over|z-b_n}~-~ \ldots ~-~
{~~a_1~~|\over|z-b_{1}} 
 \label{cftr}
\end{equation}
has the same support points as $\m_n$. In the following we are interested
into measures for which this "reversed" measure is "nearly" a uniform
distribution on its support points. More precisely, we ask
for all probability measures on $(-\infty,\infty )$ with the property 
 \begin{equation}~~~~~~~
 \left\{ \begin{array}{l}
 - \hbox{ If $n=2m-1$ is odd }(n\in \N), \hbox{ then }
\m_n^R \hbox{ has equal masses at all $n$  support }\\
 \hbox{~~~points} \\
 - \hbox{ If }n=2m \hbox{ is even }(n\in \N), \hbox{ then }
\m_n^R \hbox{ has equal masses at }n\hbox{ support points}\\
 \hbox{~~~and positive mass at a point }x_0~ (\hbox{independent of }n). \\
\end{array}
\right.
\label{ca1}
\end{equation}
The following theorem shows that there is exactly one probability
measure with the property (\ref{ca1})  (up to
a linear transformation), namely the measure with density
 proportional to the weight function of the generalized
Hermite polynomials.

\bigskip

{\bf Theorem 2.1.}~~{\it The generalized Hermite polynomials 
$\hg _n(x)$ orthogonal with respect to the
measure $d\xi_{\g}(x) = |x|^\g\exp(-x^2)dx$
 $(\g >-1)$
can be characterized as the unique (up to a linear transformation)
orthogonal polynomials on $(-\infty,\infty)$ whose corresponding
probability measure satisfies (\ref{ca1}) for all $n\in \N$. \\
Moreover,  for all $m\in \N$ the weight of 
$\xi_{\g,2m}^R$ at the point $0$ is $\g+1$ times
bigger than the (equal) weight of $\xi_{\g,2m}^R$ at 
the remaining $2m$ support points.}

\bigskip

{\bf Proof.}~~In a first step we show that the probability 
measure $\xg$ with density proportional to the function 
$|x|^\g\exp(-x^2)$
satifies (\ref{ca1}). To this end let $\hg_n (x)$ denote
the $n$th generalized Hermite polynomial and $\kg _n(x)$ its monic
form. From  [6, p. 157] we have the recursive relation
($\kg_{-1}(x)=0$, $\kg_0(x)=1$)
\begin{equation}
\kg_{n+1}(x)~=~x\kg_n(x)~-~\ah_n\kg_{n-1}(x) 
\label{rec}
\end{equation}
where
\begin{equation}
\ah_n~=~\left\{ \begin{array}{l}
\displaystyle{{n\over 2}}  \hbox{~~~~~~~~~~if }n\hbox{ is even} \\
\\
\displaystyle{{n+\g \over 2}} \hbox{~~~~ if }n\hbox{ is odd~.}\\
\end{array}
\right.
\label{an}
\end{equation}
Consequently the quantities $b_j$ in the continued fraction
expansion of the Stieltjes transform of $\xg$ in (\ref{cfg})
satisfy $b_j=0$ while the "numerators" $a_j$ are given by
(\ref{an}). Thus we obtain from (\ref{cftr})
\begin{equation}
\Phi^R_n(z)~= ~\int_{-\infty}^\infty {d\xgn^R (x)\over z - x} ~=~
{~1~|\over |~z~}~-~{\ah_n~|\over|~z~}~-~ \ldots ~-~
{\ah_1~|\over|~z~} ~~.
 \label{rev}
\end{equation}
Now let $n=2m$ ($m\in \N_0$), then it follows from
$\kg_j(z) = 2^{-j} \hg_j(z)$, (\ref{rec}), (\ref{an}) and  
formula (2.47) in [6, p. 157] that 
\begin{eqnarray}
{d\over dz}\kg_{2m+1} &~=~& (2m+1)\kg _{2m}(z)+{m\g \over z} \kg_{2m-1}
(z) \nonumber \\
&& \\
&~=~& (2m+1+\g)\kg_{2m}(z) -{\g \over z}\kg_{2m+1}(z) ~\nonumber
\label{abbll}
\end{eqnarray}
and formula (\ref{rev}) yields that the support of $~\xi_{\g,2m}^R$
is given by the zeros of $\kg_{2m+1}(z)$. Whenever
$\kg_{2m+1}(z_0)=0$ and $z_0\neq 0$  we  have from (\ref{abbll})
$$
\xi_{\g,2m}^R(z_0)~=~{\kg_{2m}(z_0) \over \left.{
d\over dz}\kg_{2m+1}(z)\right|_{z=z_0}}
~=~{1\over 2m+1+\g}~.
$$
Consequently, $\xi_{\g,2m}^R$  has equal masses ${1\over 2m+1+\g}$
at the zeros of $\kg_{2m+1}(z)$ which are different from $0$ 
and mass ${\g+1\over 2m+1+\g}$ at the point $0$.\\
In the case  $n=2m-1$ we obtain by a similar reasoning that
$$
{d\over dz}\kg_{2m}(z)~=~2m\kg_{2m-1}(z)
$$
which implies for all $z_0\in \hbox{supp}(\xi_{\g,2m-1}^R)
=\{z| ~\kg _{2m}(z) =0\}$
$$
\xi_{\g,2m-1}^R(z_0)~=~{\kg_{2m-1}(z_0) \over \left.{
d\over dz}\kg_{2m}(z)\right|_{z=z_0}}
~=~{1\over 2m}~.
$$
This proves that the probability measure with density proportional
to $|x|^\g\exp(-x^2)$ satisfies the reversing property (\ref{ca1}).\\
In a second step we now show that there is no other probability measure
with this property. If $n=2m-1$, then it follows from 
 [7, p. 16] that the property (\ref{ca1}) implies
for $m\geq 1$
\begin{eqnarray}
(2m-1)\sum_{i=1}^{2m} b_i &~=~&2m \sum_{i=1}^{2m-1}b_i \nonumber\\
&\\
(2m-2)\left[ \sum_{1\leq i<j\leq 2m} b_ib_j -\sum_{i=1}^{2m-1}a_i\right]
&~=~& 2m\left[\sum_{1\leq i<j\leq 2m-1} b_ib_j -
\sum_{i=1}^{2m-2}a_i\right]~.\nonumber
\label{eq1}
\end{eqnarray}
If $n=2m$ we denote by $x_0,x_1,\ldots ,x_{2m}$ the support
points of a reversed measure $\mu_{2m}^R$ satisfying (\ref{ca1})
(note that $x_1,\ldots ,x_{2m}$ depend on $2m$
and $x_0$ is independent of $2m$)
and obtain for the Stieltjes transform for some $\g >-1$
\begin{eqnarray*}
\Phi^R_{2m}(z)&~=~&{1\over 2m+1+\g}\sum_{j=1}^{2m}{1\over z-x_j}
~+~{\g+1\over 2m+1+\g}{1\over z-x_0} \\
&& \\
&~=~&{\prod\limits_{j=0}^{2m}(z-x_j)^{-1}\over 2m+1+\g}\Biggl[(2m+1+\g)z^{2m}
-\Bigl\{ (2m+\g)\sum_{i=0}^{2m}x_i -\g x_0 \Bigr\}z^{2m-1} \Biggr. \\
& & ~~~~~~~~~~~~~\Biggl.
+\Bigl\{(2m-1+\g)\sum_{0\leq i<j\leq 2m} x_ix_j -\g x_0 \sum_{i=1}^{2m}
x_i \Bigr\}z^{2m-2} -~\ldots \Biggr]~. 
\end{eqnarray*}
On the other hand the coefficient of $z^{2m-1}$ and $z^{2m-2}$
of the polynomial in the numerator of $\Phi^R_{2m}(z)$
in (\ref{cftr}) can be written as (see [11, p. 7]) 
$$
-\sum_{i=1}^{2m}
b_i~, ~~~ \sum_{1\leq i<j\leq 2m} b_ib_j - \sum_{i=1}^{2m-1}a_i
$$
respectively, while it is easy to see that
$$
\sum_{i=0}^{2m}x_i~=~\sum_{i=1}^{2m+1}b_i~,~~~
\sum_{0\leq i<j\leq 2m} x_ix_j~=~\sum_{1\leq i<j\leq 2m+1}b_ib_j
-\sum_{i=1}^{2m}a_i~.
$$
Thus the property (\ref{ca1}) (in the case $n=2m$, $m\geq 0$) 
implies  the equations
\begin{eqnarray*}
(2m+\g)\sum_{i=1}^{2m+1}b_i -\g x_0 ~=~(2m+1+\g)\sum_{i=1}^{2m}  
b_i && \\
&& \\
(2m-1+\g) \left[\sum_{1\leq i<j\leq 2m+1}b_ib_j-\sum_{i=1}^{2m}a_i\right]
-\g x_0\left[\sum_{i=1}^{2m+1}b_i-x_0\right]~&=&~ \\
 (2m+1+\g)
\left[\sum_{1\leq i<j\leq 2m} b_ib_j \right .&-&
\left .\sum_{i=1}^{2m-1}a_i\right]
\end{eqnarray*}
which can easily be rewritten as
\begin{eqnarray}
b_{2m+1} &~=~& \displaystyle{1\over 2m+\g} 
\left[\sum_{i=1}^{2m} b_i +\g x_0\right] 
\nonumber \\
&&\\
(2m-1+\g)a_{2m}&~=~& 2\sum_{i=1}^{2m-1}a_i - \g x_0 \left[\sum_{i=1}^{2m+1} 
b_i-x_0\right] - 2 \sum_{1\leq i< j \leq 2m} b_ib_j \nonumber \\ 
&&+ ~~~(2m-1+\g)b_{2m+1}\sum_{i=1}^{2m}b_i \nonumber 
\label{eq2}
\end{eqnarray}
($m\geq 0$). From the first equations in (\ref{eq1}) and (\ref{eq2})
it follows that $b_j=x_0$ for all $j\in \N$ while the second 
equations reduce to
\begin{eqnarray}
a_{2m-1} &~=~& {2\over 2m-2} \sum_{j=1}^{2m-2} a_j ~~~~~~~~~~~~~~~m \geq 2 
\nonumber \\
&&  \\
a_{2m} &~=~& {2\over 2m-1+\g} \sum_{j=1}^{2m-1} a_j~~~~~~~~~~m \geq 1~.
\nonumber 
\label{eq3}
\end{eqnarray}
An induction argument now shows that all solutions of (\ref{eq3}) are
of the form
\begin{equation}
a_j~=~\left\{\begin{array}{l}
 \displaystyle{j\over 2}  c^{-2}~
 \hbox{~~~~~~~~~~~~if } j \hbox{ is even} \\
 \\
\displaystyle{j+\g \over 2}  c^{-2}
 \hbox{~~~~~~~ if } j \hbox{ is odd} \\
\end{array}~~~~~(j \in \N,~ c>0)\right.
\label{canmom}
\end{equation}
and by the discussion in the first part of the proof
the measure $\m$ corresponding to this sequence
has the density proportional to $|c(x-x_0)|^\g \exp(-|c(x-x_0)|^2)$.
This completes the proof of Theorem 2.1. \hfill \bull

 \bigskip
For $\g=0$ Theorem 2.1 reduces to the characterization of the
Hermite polynomials given in [7]. In this case the measure 
$\xgn ^R$ puts equal masses at its support points or all $n\in
\N$. The following result provides a slightly different characterization
of the generalized Hermite polynomials in the class of all
symmetric polynomials.

\bigskip

{\bf Theorem 2.2.}~~{\it The generalized Hermite polynomials
$\hg _n(x)$ orthogonal with respect to the measure
$d\xi _\g(x) =|x|^\g\exp (-x^2)dx$  $(\g >-1)$ 
can be characterized as the unique (up to a
scaling factor) symmetric orthogonal polynomials whose
corresponding probability measure $\m$  satisfies for all $m\in \N$
\begin{equation} 
\m_{2m}^R \hbox{ has equal masses at all support points which are different
from zero}
\label{ca2}
\end{equation}
Moreover, for all $m \in \N$ the weight of 
$\xi_{\g,2m}^R$ at the point $0$ is $\g+1$ times
bigger than the (equal) weight of $\xi_{\g,2m}^R$ at 
the remaining $2m$ support points.}

\bigskip

{\bf Proof.}~~The proof is similar to that of Theorem 2.1 and we
only sketch the main steps. By the discussion in
the first part of the proof of the previous theorem it follows
that the measure $\xg$ with density proportional
to $|x|^\g \exp(-x^2)$ satisfies (\ref{ca2}).
Assume now that $\m$ is a symmetric measure
such that (\ref{ca2}) holds for all $m\in \N$. By the
symmetry of $\m$ we have $b_i=0$ ($i\in \N$) and for some $\g > -1$
\begin{eqnarray*}
&&\Phi^R_{2m}(z)~=~\int_{-\infty}^\infty {d\mu^R_{2m} (x)\over z - x} ~=~
{~1~|\over |~z~}~-~{a_{2m}|\over|~z~}~-~ \ldots ~-~
{a_1|\over|~z} ~=~z^{-1}\prod_{j=1}^{2m}(z-x_j)^{-1}\\
&~\times&\left[ z^{2m}-
z^{2m-2} {2m-1+\g\over 2m+1+\g}\sum_{1\leq i<j\leq 2m} x_ix_j
+ z^{2m-4} {2m-3+\g\over 2m+1+\g} \sum_{1\leq i<j<k<l\leq 2m}
x_ix_jx_kx_l ~\ldots \right] 
\end{eqnarray*}
where $x_0=0$, $x_1,\ldots x_{2m}$ denote the
$2m+1$ support points of the (symmetrc) measure 
 $\m _{2m}^R$. Comparing coefficients
of the two representations for the polynomial in the denominator
of this  continued fraction (similary as in the second 
part of the proof of Theorem 2.1)
now yields the equations
\begin{eqnarray}
a_{2m} &~=~&{2\over 2m-1+\g}\sum_{i=1}^{2m-1}a_i 
~~~~~~~~~~~~~~~~~~~~~~~~m\geq 1 \nonumber\\
&&\\
a_{2m}\sum_{i=1}^{2m-2} a_i &~=~&{4\over 2m-3+\g}\sum_{1\leq
i <j \leq 2m-2} a_ia_{j+1} ~~~~~~~~~~~m\geq 2. \nonumber
\label{eq4}
\end{eqnarray}
A tedious calculation shows that (\ref{eq4}) is equivalent
to ($a_0=0$)
\begin{eqnarray*}
a_{2m} &~=~&{2\over 2m-1+\g}a_{2m-1}+a_{2m-2}
~~~~~~~~~m\geq 1 \\
&&\\
a_{2m}&~=~&{4\over 2m-1+\g}a_{2m-1}+ a_{2m-4}~~~~~~~~~m\geq 2
\end{eqnarray*}
and an induction argument shows that the (unique) solution
of this system is given by (\ref{canmom}). This
 completes the proof of the theorem. \hfill \bull

\bigskip

In the remaining part of this Section we will concentrate
 on symmetric distributions on the interval $[-1,1]$ 
 for which similar characterizations can be derived. For the
 sake of brevity we only consider a characterization
of the type (\ref{ca2}). Let $\m$ denote a symmetric
probability measure on the interval $[-1,1]$ with Stieltjes 
transform
\begin{equation}
\Phi (z)~= ~\int_{-1}^1 {d\mu (x)\over z - x} ~=~
{~1~|\over |~z~}~-~{~p_2~|\over|~z~}~-~{q_2p_4|\over|~~z~} 
~-~{q_4p_6|\over|~~z~}~-~ \ldots
\label{int}
\end{equation}
where $p_{2i} \in [0,1]$, $q_{2i}=1-p_{2i}$ ($i\in \N$),
$q_0=1$ and $(q_{2i-2}p_{2i})_{i\in \N}$ is the minimal chain
sequence in the recursive relation of the monic
orthogonal polynomials with respect to the measure $d\m (x)$.
If (\ref{int}) holds for a symmetric
probability measure $\m$ (on the
interval $[-1,1]$) we write that "$\m$ corresponds to the
sequence $(p_2,p_4,p_6, \ldots )$". 
If $\m _n$ denotes the (finite) measure corresponding
to the "terminated" continued fraction
\begin{equation}
\Phi_n (z)~= ~\int_{-1}^1 {d\mu _n(x)\over z - x} ~=~
{~1~|\over |~z~}~-~{~p_2~|\over|~z~}~-~{q_2p_4|\over|~~z~}~-~\ldots 
~-~{q_{2n-2}p_{2n}|\over|~~~~z~~~~} 
\label{term}
\end{equation}
we write that "$\m_n$ corresponds to the  sequence
$(p_2,\ldots ,p_{2n},0)$" and denote by
$\m_n^R$ the measure
corresponding to the sequence $(p_{2n},\ldots ,p_2,0)$.
The Stieltjes transform of $\m_n^R$ is given by
\begin{equation}
\Phi _n^R (z)~= ~\int_{-1}^1 {d\mu _n^R (x)\over z - x} ~=~
{~1~|\over |~z~}~-~{~p_{2n}|\over|~~z~}~-~{q_{2n}p_{2n-2}|\over|~~~~z~~~~}
~ -~\ldots ~-~{q_{4}p_{2}|\over|~z~~} 
\label{termrev}
\end{equation}
and it is known that $ \m_n$ and $\m_n^R$ have the same support
points (see [7]).\\
Throughout this paper we define  $\cag_n(x)$ as the 
"generalized ultraspherical polynomials" orthogonal
with respect to the measure $|x|^\g (1-x^2)^{\a-1/2}dx$
(without specifying any normalization). 
The following result is the analogue of Theorem 2.2
on the interval $[-1,1]$ and can be proved by similar 
arguments as given in the proof of  Theorem 2.2.

\bigskip

{\bf Theorem 2.3.}~~{\it The generalized ultraspherical 
polynomials $\cag _n(x)$ orthogonal with respect to the
measure $d\xi_{\a,\g}^*(x)=|x|^\g(1-x^2)^{\a-1/2} dx$ 
$(\a >-1/2$, $\g >-1$) can be characterized as the unique
symmetric orthogonal polynomials whose corresponding
probability measure $\m$ satisfies for all $m\in \N$
\begin{equation}
\m _{2m}^R \hbox{ has equal masses at all support points
different from } 0.
\label{ca3}
\end{equation}
Moreover, for all $m\in \N$ the weight of 
$\xi_{\a,\g,2m}^{R^*}$ at the point $0$ is $\g+1$ times
bigger than the (equal) weight of $\xi_{\a,\g,2m}^{R^*}$ at 
the remaining $2m$ support points.}

\bigskip

{\bf Proof.}~~We will only show that the probability measure
$\xi_{\a,\g}^*$ with density proportional to the function 
$|x|^\g(1-x^2)^{\a-1/2}$ has the property (\ref{ca3}). The
result that there are no other measures with this property
can be shown exactly in the same way as the 
corresponding statements in Theorem 2.1 and 2.2. \\
The monic polynomials $\sh^{(\a,\g)} _n(x)$ orthogonal with
respect to the measure  $|x|^\g(1-x^2)^{\a-1/2}dx$ can be
obtained from  [6, p. 156] as $\sh^{(\a,\g)} _{-1}(x)=0$,
$\sh^{(\a,\g)} _0(x)=1$,
\begin{equation}
\sh^{(\a,\g)} _{n+1} (x) ~=~ x\sh^{(\a,\g)} _n(x) ~-~ 
\g_{n+1}^{(\a,\g)} \sh^{(\a,\g)}_{n-1} (x)
\label{defin}
\end{equation}
where 
\begin{eqnarray}
\g _{2m}^{(\a,\g)} &~~=~\displaystyle{(2m-1+\g)(2m-2+2\a+\g) \over
(4m-4+2\a+\g)(4m-2+2\a+\g)} \nonumber\\
&& \\
\g _{2m+1}^{(\a,\g)} &=~\displaystyle{2m(2m+2\a-1) \over
(4m-2+2\a+\g)(4m+2\a+\g)} ~\nonumber
\label{rec1}
\end{eqnarray}
Consequently $\g _{n+1}^{(\a,\g)} =q_{2n-2}p_{2n}$  is a chain sequence
with ($q_0=1$)
\begin{equation}
p_{2j}~=~\left\{\begin{array}{l}
\displaystyle{j\over 2\a+\g+2j}~~~~~~\hbox{if } j \hbox{ is even}\\
\\
\displaystyle{j+\g \over 2\a+\g+2j}~~~~~~\hbox{if } j \hbox{ is odd}~.
\end{array}\right.
\label{pn}
\end{equation}
Thus we have form (\ref{termrev}), (\ref{rec1}), (\ref{pn}) and
straightforward calculations that 
\begin{eqnarray}
\Phi^{R^*}_{2m}(z)~&=&~\int_{-1}^1 {d\xgam^{R^*} (x)\over z - x} ~=~
{~1~|\over |~z~}~-~{p_{4m}|\over|~~z~}~-~{q_{4m}p_{4m-2}|\over
|~~~~~z~~~~ }~-~ \ldots ~-~{q_4p_2|\over |~z~~}\nn \\
&& \\
&=&~{~1~|\over |~z~}~-~{p_{4m}|\over|~z~}~-~{\g_{2m}^{(\a+1,\g)}|\over
|~~~~z~~~ }~-~\ldots ~-~{\g_2^{(\a+1,\g)}|\over |~~~~z~~~} \nn\\
&&\nn\\
~&=&~\displaystyle{\sh _{2m}^{(\a+1,\g)}(z) \over
z\sh _{2m}^{(\a+1,\g)}(z) - p_{4m} \sh _{2m-1}^{(\a+1,\g)}(z)}~~.\nn
\label{lab}
\end{eqnarray}
Using the relations
\begin{eqnarray*}
\sh _{2m}^{(\a+1,\g)}(z) ~&=&~\displaystyle{\G (m+1) \G (m+\a +1 +{\g\over 2} )
\over \G (2m +\a +1+{\g \over 2})} ~P_{m}^{(\a+1/2,(\g-1)/2)}(2z^2-1) \\
&& \\ 
\sh _{2m-1}^{(\a+1,\g)}(z) ~&=&~\displaystyle{\G (m) \G (m+\a +1 +{\g\over 2} )
\over \G (2m +\a +{\g \over 2})} ~zP_{m-1}^{(\a+1/2,(\g+1)/2)}(2z^2-1) 
\end{eqnarray*}
(see  [6, p.156]), (\ref{pn}), formula (22.7.20)
in Abramowitz and Stegun [1] it now follows that
\begin{eqnarray*}
&&z \sh _{2m}^{(\a+1,\g)}(z)~-~p_{4m}\sh_{2m-1}^{(\a+1,\g)}(z)  \\
&&\\
&=&~ z\displaystyle{\G (m+1) \G (m+\a +1 +{\g\over 2} ) \over \G (2m+\a+1+
{\g\over 2})} \left[ P_{m}^{(\a+1/2,(\g-1)/2)}(2z^2-1) -
P_{m-1}^{(\a+1/2,(\g+1)/2)}(2z^2-1) \right]\\
&&\\
&=&~ z\displaystyle{\G (m+1) \G (m+\a +1 +{\g\over 2} ) \over \G (2m+\a+1+
{\g\over 2})}  P_{m}^{(\a-1/2,(\g+1)/2)}(2z^2-1)
\end{eqnarray*}
and the Stieltjes transform of $\xgam ^{R^*}$ is obtained as
\begin{equation}
\Phi^{R^*}_{2m}(z)~=~ {P_{m}^{(\a+1/2,(\g-1)/2)}(2z^2-1)
\over z P_{m}^{(\a-1/2,(\g+1)/2)}(2z^2-1)} ~.
\label{stie}
\end{equation}
Consequently the support points of $\xgam ^{R^*}$ are given by the
polynomial in the denominator of (\ref{stie}) and we have
from Szeg\"o [13, p. 63] for all
nonvanishing support points $z_0$ of $\xgam  ^{R^*}$
\begin{eqnarray*}
&& {d\over dz} \left. \left(zP_m^{(\a-1/2,(\g+1)/2)}(2z^2-1) \right)
\right|_{z=z_0}~\\
&&\\
&=&~P_m^{(\a-1/2,(\g+1)/2)}(2z_0^2-1)
+(2m+2\a+\g+2) z_0^2P_{m-1}^{(\a+1/2,(\g+3)/2)}(2z_0^2-1) \\
&&\\
&=&~\displaystyle{m+\a+{\g\over 2}+1\over
m+{\a+1\over 2}+{\g\over 4}}
\left[(m+{\g+1\over 2}) P_m^{(\a+1/2,(\g-1)/2)}(2z_0^2-1)
\right. \\
&-&~\left.(m+{\g+1\over 2}) P_m^{(\a-1/2,(\g+1)/2)}(2z_0^2-1)
+m P_m^{(\a+1/2,(\g+1)/2)}(2z_0^2-1) \right]\\
&&\\
&=&~(2m+\g+1)P_m^{(\a+1/2,(\g-1)/2)}(2z_0^2-1)~.
\end{eqnarray*}
Here the second equality follows from the first by the formulas
(22.7.16, 22.7.20) in [1] while
the last equality follows from the second by formula
(22.7.19, 22.7.20) (in the same reference) and $P_m^{(\a-1/2,(\g+1)/2)}
(2z_0^2-1)=0$. Thus we obtain from (\ref{stie}) for every nonvanishing
support point $z_0$ of $\xgam ^{R^*}$
$$
\xgam ^{R^*}(z_0)~=~{P_m^{(\a+1/2,(\g-1)/2)}(2z_0^2-1)\over 
{d\over dz}\left.\left( P_m^{(\a-1/2,(\g+1)/2)}
(2z^2-1)\right)\right|_{z=z_0}}~=~{1\over 2m+\g+1}
$$
which proves the assertion of the Theorem.\hfill \bull

\bigskip

\begin{center}
  {\sc 3. Generalized sieved ultraspherical polynomials}\\
\end{center}
\def\theequation{3.\arabic{equation}}
\setcounter{equation}{0}

In Theorem 2.3 of Section 2  we considered the $1$th, $3$th,
$5$th, $7$th $~\ldots ~$ convergent of the
continued fraction expansion for the Stieltjes transform
of a given probability measure on the interval $[-1,1]$,
reversed the sequences of corresponding $p_{2i}$'s and
determined all measures with the property (\ref{ca3}).
A natural extension of this proceeding is to terminate the
continued fraction at the positions $k$, $3k$, $5k$, $~\ldots ~$ 
(for a given $k\in \N$) and  investigate
if there exist similar characterizations.
We will show in this Section that these problems are related to 
the so called "sieved random walk polynomials".\\
 Following the work of [4], [5]
and [10] we use a set of "random walk polynomials"
\begin{equation}
\left\{\begin{array}{l}
xR_n(x)~=~A_nR_{n+1}(x)+B_nR_{n-1}(x)~~~~~~n\geq 0 \\
\\
R_{-1}(x)=0,~R_0(x)=1 
\end{array}\right.
\label{RW}
\end{equation}
$(A_n,B_n>0$, $A_n+B_n=1$, $n\geq 0$) in order to define "sieved random walk
polynomials" of the first kind" by
\begin{equation}
\left\{\begin{array}{l}
xr_n(x)~=~b_{n-1}r_{n+1}(x)+a_{n-1}r_{n-1}(x)~~~~~~~n\geq 1 \\
\\
r_{1}(x)=x,~r_0(x)=1 
\end{array}\right.
\label{SRW1}
\end{equation}
and "sieved random walk polynomials of the second kind" by
\begin{equation}
\left\{\begin{array}{l}
xs_n(x)~=~a_{n}s_{n+1}(x)+b_{n}s_{n-1}(x)~~~~~~n\geq 1 \\
\\
s_{1}(x)=2x,~s_0(x)=1 
\end{array}\right.
\label{SRW2}
\end{equation}
where 
\begin{equation}
a_n=b_n={1\over 2}~~\hbox{ if } k\not ~\mid n+1, ~a_{nk-1}=A_{n-1},~B_{nk-1}=B_{n-1}
\label{coef}
\end{equation}
and $k\geq 2$ is a fixed integer.
 In the following we are interested
 into the polynomials $r_n(x)$ and $s_n(x)$ when one uses the
 generalized ultraspherical polynomials $C_n^{(\a+1,\g)}(x)$ 
 defined in Section 2 as
 random walk polynomials. From [6, p.156]
 we see that in this case
 \begin{eqnarray}
 A_{2m}~&=&~1-B_{2m}~=~\displaystyle{2m+1+\g \over 4m+2\a+\g+2} \nn \\
 && \\
 A_{2m-1} ~&=&~1-B_{2m-1} ~=~\displaystyle{2m \over 4m+2\a+\g}\nn
 \label{Am}
 \end{eqnarray}
 (note that the $A_n$'s correspond to the $p_{2n}$'s in 
(\ref{pn}), i.e. $A_{n-1}=p_{2n}$).
Throughout this paper we will denote the polynomials obtained from
(\ref{RW}), (\ref{SRW1}), (\ref{SRW2}), (\ref{coef}) and
(\ref{Am}) as
"\gsu polynomials of the first kind" $\cag _n (x,k)$
and "\gsu  polynomials of the second kind" $\bag _n(x,k)$.
A straightforward calculation yields
 the recursive relations $\cag _0(x,k)=1$,
$\cag _1(x,k)=x$,
\begin{equation}
\left\{ \begin{array}{l}
(2\a+\g+2n)x\cag _{nk}(x,k) =(2\a+\g+n) \cag_{nk+1} (x,k)+n
\cag_{nk-1} (x,k)~
\\~~~~~~~~~~~~~~~~~~~~~~~~~~~~~~~~~~~~~~~~~~~
~~~~~~~~~~~~~~~~~~~~~~~~~~~~~~~~~~~~\hbox{ if }n\hbox{ is even} \\
(2\a+\g+2n)x\cag _{nk}(x,k) =(2\a+n) \cag_{nk+1} (x,k)+(n+\g)
\cag_{nk-1} (x,k)\\~~~~~~~~~~~~~~~~~~~~~~~~~~~~~~~~~~~~~~~~~~~
~~~~~~~~~~~~~~~~~~~~~~~~~~~~~~~~~~~~~\hbox{if }n\hbox{ is odd} \\
2x\cag _{j}(x,k) ~=~\cag_{j+1} (x,k)+
\cag_{j-1} (x,k)~~~~~~\hbox{if }j\neq nk
\end{array}
\right.
\label{su1}
\end{equation}

\medskip

and $\bag _0(x,k)=1$, $\bag _1(x,k)=2x$,

\medskip

\begin{equation}
\left\{ \begin{array}{l}
(2\a+\g+2n)x\bag _{nk-1}(x,k) =n \bag_{nk} (x,k)+(2\a+\g+n)
\bag_{nk-2} (x,k)\\~~~~~~~~~~~~~~~~~~~~~~~~~~~~~~~~~~~~~~~~~~~
~~~~~~~~~~~~~~~~~~~~~~~~~~~~~~~~~~~~\hbox{if }n\hbox{ is even} \\
(2\a+\g+2n)x\bag _{nk-1}(x,k) =(n+\g) \bag_{nk} (x,k)+(n+2\a)
\bag_{nk-2} (x,k)\\~~~~~~~~~~~~~~~~~~~~~~~~~~~~~~~~~~~~~~~~~~~
~~~~~~~~~~~~~~~~~~~~~~~~~~~~~~~~~~~~\hbox{if }n\hbox{ is odd} \\
2x\bag _{j}(x,k) ~=~\bag_{j+1} (x,k)+
\bag_{j-1} (x,k)~~~~~~~~~\hbox{if }j+1\neq nk~.
\end{array}
\right.
\label{su2}
\end{equation}
Note that the case $\g=0$ gives the sieved ultraspherical polynomials
which were obtained in [2]
as a limit from the $q-$ultraspherical polynomials of [12].
The corresponding measure of orthogonality can be
obtained from [9, p. 561] and [10, p. 96] as
\begin{equation}
w_1(x,\a,\g) ~=~(1-x^2)^{\a-1/2}|U_{k-1}(x)|^{2\a}|T_k(x)|^\g~~~~~~(x \in [-1,1])
\label{w1}
\end{equation}
for the \gsu polynomials of the first kind and as
\begin{equation}
w_2(x,\a,\g) ~=~(1-x^2)^{\a+1/2}|U_{k-1}(x)|^{2\a}|T_k(x)|^\g~~~~~~(x \in [-1,1])
\label{w2}
\end{equation}
for the polynomials of the second kind.
The following Theorem characterizes the \gsu polynomials of the
first kind by a similar "reversing property" as stated in Section 2.

\bigskip

{\bf Theorem 3.1.}~~{\it The \gsu polynomials of the first kind
orthogonal with respect to the measure $d\xi_{\a,\g}(x)
=w_1(x,\a,\g)dx$ defined by (\ref{su1}) 
can be characterized as the unique sieved random walk polynomials
of the first kind whose corresponding probability measure $\m$
satisfies for all $m\in \N_0$
\begin{equation}
~~~~~~~~~~~
\left\{ \begin{array}{l} 
- \m ^R_{k(2m+1)-1} \hbox{has equal masses at the zeros of }
T_k(x) \\
- \m ^R_{k(2m+1)-1}  \hbox{ has equal (but not necessary the same) masses
 at the remaining }\\
~~2mk  \hbox{ support points}
 \end{array}\right.
 \label{ca4}
 \end{equation}
Moreover, for all $m \in  \N$ the masses of 
$\xi_{\a,\g,k(2m+1)-1}^{R}$ at the zeros of $T_k(z)$
are  $\g+1$ times
bigger than the (equal) masses of $\xi_{\a,\g,k(2m+1)-1}^{R}$ at 
the remaining $2mk$ support points.}

 \bigskip

{\bf Proof.}~~Consider a set of random walk polynomials  defined
by (\ref{RW}) and the corresponding set of sieved polynomials
of the first kind in (\ref{SRW1}) orthogonal with respect
to the measure $\m$. Observing (\ref{SRW1}) and (\ref{coef})
 it is straightforward to see that the  minimal chain
 sequence $(q_{2i-2}p_{2i})_{i\in \N}$
$(q_0=1)$ corresponding to the monic orthogonal polynomials satisfies
\begin{equation}
p_{2i}~=~{1\over 2}~~~~~~~~i \neq nk
\label{lab3}
\end{equation}
and we obtain from (\ref{int}) for the Stieltjes transform of $\m$
\begin{eqnarray*}
&&\Phi (z) ~=~\int_{-1}^1 {d\m (x) \over z-x}~=~
{~1~|\over |~z~}~-~{~p_2~|\over|~z~}~-~{q_2p_4|\over|~~z~} 
~-~{q_4p_6|\over|~~z~}~-~\ldots \\
&&\\
=&&{~1~|\over |~z~}-{~{1\over 2}~|\over |~z~}-
\underbrace{{{1\over 4}~|\over|~z~}-
\ldots  -{{1\over 4}~|\over|~z~}}_{k-2} - {{1\over 2}p_{2k}|\over |~~z~}
-{{1\over 2}q_{2k}|\over |~~z~}
-\underbrace{{{1\over 4}~|\over|~z~}-\ldots  -{{1\over 4}~|\over|~z~}}_{k-2}
-{{1\over 2}p_{4k}|\over|~~z~} -{{1\over 2}q_{4k}|\over|~~z~}
-~\ldots 
\end{eqnarray*}
A contraction such that the resulting continued fraction
attains successively the values of the $(k-1)$th, $(2k-1)$th,
$(3k-1)$th $\ldots $ convergent  yields (see [11], p. 12)
\begin{equation}
\Phi (z) ~=~\left( z-{1\over 2}H(z) \right)^{-1}
\label{b1}
\end{equation}
where
\begin{eqnarray}
H(z) &=&2{U_{k-2}(z)\over U_{k-1}(z)}+
{~~~~~~~({1\over 2})^{2k-3} p_{2k}~~~~~|\over |({1\over 2})^{2k-2}
U_{k-1}(z)T_k(z)}- 
{({1\over 2})^{3k-3} q_{2k}p_{4k}U_{k-1}(z)|\over |~~~~~({1\over 2})^{k-1}
T_k(z)~~~~~~}-
{~({1\over 2})^{2k-2} q_{4k}p_{6k}|\over |~({1\over 2})^{k-1}
T_k(z)~~}~\ldots \nonumber \\
&&\\
&=&~{2\over U_{k-1}(z)} \left[U_{k-2}(z)
~+~ {~p_{2k}~|\over |T_k(z)}~-~{q_{2k}p_{4k}|\over|T_k(z)}~-
~{q_{4k}p_{6k}|\over|
T_k(z)} ~-~ \ldots ~~\right] ~.\nonumber
\label{b2}
\end{eqnarray}
Thus we obtain from (\ref{b1})
\begin{eqnarray}
\Phi (z) ~&=&~\int_{-1}^1 {d\m (x) \over z-x}~=~
 U_{k-1}(z) \left[~{~~1~~|\over |T_k(z)} -
 {~p_{2k}~|\over |T_k(z)}-{q_{2k}p_{4k}|\over|T_k(z)}-{q_{4k}p_{6k}|\over|
T_k(z)} -~ \ldots ~\right] \nonumber\\
&& \\
&=&~ U_{k-1}(z) 
\int_{-1}^1 {d\ms (x) \over T_k(z)-x}\nonumber
\label{b3}
\end{eqnarray}
where $\ms$ is the probability measure defined by 
the recursive relation of the monic orthogonal
polynomials $R'_{-1}(x)=0,~R'_0(x)=1$ 
\begin{equation}
\left\{\begin{array}{l}
R'_{n+1}(x)~=~xR'_n(x)- q_{2(n-1)k}p_{2nk}R'_{n-1}(x)\\ 
\\
~~~~~~~~~~~=xR'_n(x)- B_{n-2}A_{n-1}R'_{n-1}(x)
\end{array}\right.
\label{b4}
\end{equation}
($q_0=1$, $A_{-1}=0$), see also [5, p. 82] or [9, p. 562]. \\
Now assume that $\xi_{\a,\g}$ is the probability 
 measure with density proportional
to the weight function ({\ref{w1}) such that the \gsu polynomials 
$\cag _n(x,k)$ of the first kind are
orthogonal with respect to the measure 
$d\xi_{\a,\g} (x)$. From (\ref{Am})
it follows that the polynomials $R_l'(x)$ coincide with
the polynomials $\sh _l^{(\a,\g)}(x)$ in (\ref{defin}).
Therefore the measure $\xi_{\a,\g} ^*$ defined
by (\ref{b4}) corresponds to the sequence (\ref{pn})
and has density proportional to $|x|^\g(1-x^2)^{\a-1/2}$.
A similar reasoning as in the derivation of 
(\ref{b2}) yields
\begin{eqnarray*}
\Phi^R _{k(2m+1)-1}(z)&~=~&\int_{-1}^1 {d\xi^{R}_{\a ,\g ,k(2m+1)-1}(x)
\over z-x}~=~U_{k-1}(z) \int_{-1}^1 {d\xgam ^{R^*}(x)
\over T_k(z)-x} \\
&&\\
&~=:~&U_{k-1}(z) \Phi_{2m}^{R^*} (T_k(z))
~=~U_{k-1}(z) {V_{2m}(T_k(z))\over W_{2m+1}(T_k(z))}~
\end{eqnarray*}
where $\xi_{\a,\g,2m}^{R^*}$ is defined in (\ref{lab})
and  $V_{2m}(z)$, $W_{2m+1}(z)$ are polynomials of degree
$2m$ and $2m+1$, respectively.
Consequently  Theorem 2.3 yields for all $z_0 \in$ supp $(
\xi ^R_{\a,\g,k(2m+1)-1})$
\begin{eqnarray*}
&&\xi ^R_{\a,\g,k(2m+1)-1}(z_0) ~=~\left .
(z-z_0)\Phi^R_{k(2m+1)-1}(z)\right|_{z=z_0}
~=~\left. (z-z_0)U_{k-1}(z_0)\Phi_{2m}^{R^*}(T_k(z))\right|_{z=z_0} \\
&&\\
&=&~\displaystyle{{U_{k-1}(z_0)V_{2m}(T_k(z_0))\over \left.
{d\over dz}W_{2m+1}(T_k(z))\right|_{z=z_0}}
~=~{1\over k} \xgam ^{R^*}  (T_k(z_0)) ~=~
\left\{ \begin{array}{l}
{1\over k} {\g+1\over 2m+1+\g} ~~\hbox{ ~if~ } T_k(z_0)=0 \\
\\
{1\over k} {1\over 2m+1+\g} ~~\hbox{ ~if ~} T_k(z_0) \neq 0 
\end{array}\right.}
\end{eqnarray*}
which shows that $\xi_{\a,\g,k(2m+1)-1}^R$ satisfies (\ref{ca4})
for all $m\in \N_0$. \\
On the other hand, if $\m $ denotes a measure
of orthogonality for a set of sieved random walk 
polynomials of the first kind satisfying (\ref{ca4}), then
similar arguments as given above show that the measure
$\ms$ defined by  (\ref{b3}) and (\ref{b4})
satisfies (\ref{ca2}), which determines $\ms$ uniquely
in the set of all symmetric probability measures on
the interval $[-1,1]$
(by Theorem 2.3). Therefore (observing (\ref{lab3})) $\m$
is  unique and the assertion of the theorem follows.\hfill \bull

\bigskip

The following result states an analogous characterization
for the \gsu polynomials of the second kind. Its
proof is performed by similar arguments as given
in the proof of Theorem 3.1 and therefore omitted.

\bigskip

{\bf Theorem 3.2.}~~{\it The \gsu polynomials
of the second kind orthogonal with respect to the
measure $d\eta_{\a,\g}(x) =w_2(x,\a,\g)dx$
defined in (\ref{su2})  can be characterized as the unique
sieved random walk polynomials of the second kind
whose corresponding probability measure $\m$ 
satisfies for all $m \in \N_0$
\begin{equation}
\left\{ \begin{array}{l}
- \m _{k(2m+2)-2}^R \hbox{~~has equal masses at the zeros of }T_k(z)\\
- \m _{k(2m+2)-2}^R \hbox{~~has equal masses at the zeros of }U_{k-1}(z)\\
- \m _{k(2m+2)-2}^R \hbox{~~has equal masses at the remaining }2mk
\hbox{ support points} 
\end{array}\right.
\label{ca5}
\end{equation}
Moreover, for all $m \in \N$ the masses of 
$\eta_{\a,\g,k(2m+2)-2}^{R}$ at the zeros of $U_{k-1}(z)$ and 
$T_k(z)$ are $2\a+1$ and  $\g+1$ times
bigger than the (equal) masses of $\eta_{\a,\g,k(2m+2)-2}^{R}$ at 
the remaining $2mk$ support points.}

\bigskip

{\bf Remark 3.3.}~~It is worthwhile to mention that there exist 
a couple of similar properties of the weight functions
(\ref{w1}), (\ref{w2}) of the \gsu polynomials of the first and
second kind. For example, if $\xga$ is the probability measure
with density proportional to (\ref{w1}), then it can
be shown that for all $m\in \N$
$$
- \xi _{\a,\g,2mk-1}^R~~\hbox{has equal masses at all } 2mk
\hbox{ support points}.
$$
Similary it follows for the measure $\eta _{\g,\a}$ with
density proportional to (\ref{w2}) that for all
$m \in \N_0$

\medskip

$-$ $\eta_{\g,\a,k(2m+1)-2}^R $ has equal masses ${2\a+1\over 2mk+
(2\a+1)(k-1)}$ at the zeros of $U_{k-1}(z)$

\smallskip

$-$ $\eta_{\g,\a,k(2m+1)-2}^R $ has equal  masses
${1\over 2mk + (2\a+1)(k-1)}$ at the remaining $2mk$
 support points. 

\medskip

 Finally we remark that we conjecture that the properties
 (\ref{ca4}) and (\ref{ca5}) characterize the
 \gsu poynomials of the first and of the second kind in the 
 class of all symmetric
 orthogonal polynomials on the interval $[-1,1]$.

\bigskip

\begin{center}
  {\sc 4. Asymptotic distribution of the zeros}\\
\end{center}
\def\theequation{4.\arabic{equation}}
\setcounter{equation}{0}

In this Section we investigate the asymptotic behaviour of the zeros 
of the \gsu polynomials and generalized Hermite polynomials 
defined in Section 3 and 2 when the degree and the parameters tend to 
infinity. To this end consider a set of random walk
polynomials defined by (\ref{RW}) and (\ref{coef}) where
the parameters $\a$ and $\g$ depend on the degree $n$, that
is
\begin{equation}
A_j^{(n)}~=~1-B_j^{(n)}~=~\left\{ \begin{array}{l}
\displaystyle{j+1+\g_n\over 2\a_n+\g_n+2j+2} 
~~~~\hbox{if }j \hbox{ is even} \\
\\
\displaystyle{j+1\over 2\a_n+\g_n+2j+2} ~~~~\hbox{if }j \hbox{ is odd} 
\end{array}\right. ~~~~j=0,\ldots , n-1 
\label{41}
\end{equation}
and the corresponding sieved random walk polynomials 
$r_l(x)$ of 
the first kind where we use the set
$\{A_0^{(n)},A_1^{(n)},\ldots ,A_{n-1}^{(n)}\}$ in
(\ref{SRW1}) if $nk+1 \leq l \leq (n+1)k$. Thus the
parameters of the \gsu polynomials of the 
first kind depend on the degree and 
we obtain $r_l(x) =\cags (x,k)$ where
\begin{equation}
\left\{ \begin{array}{l}
\as _l ~=~ \a_{\lfloor (l-1)/k \rfloor} \\
\gs _l ~=~ \g_{\lfloor (l-1)/k \rfloor} 
\end{array}\right.
\label{astern}
\end{equation}
and $\lfloor x \rfloor$ denotes the largest integer less
or equal than $x$. For $y \in [-1,1]$ let
\begin{equation}
N_l^{(\as_l,\gs_l)}(y)~:=~\#\{x\leq y~|~\cags (x,k)=0~\}
\label{41a}
\end{equation}
denote the number of zeros of $\cags (x,k)$
less or equal than $y$. In order to
discuss the asymptotic behaviour of $N_l^{(\as_l,\gs_l)}$
when $\a_l \to \infty$, $\g_l \to \infty$ we 
need the following auxiliary result.

\bigskip

{\bf Lemma 4.1.}~~{\it Let $\rho$ denote 
a symmetric probability measure 
on the interval $[-1,1]$ with Stieltjes transform
$$
\Phi ^* (z)~= ~\int_{-1}^1 {d\rho (x)\over z - x} ~=~
{~1~|\over |~z~}~-~{~p_2|\over|~z~}~-~{q_2p_4|\over|~~z~} 
~-~{q_4p_6|\over|~~z~} ~-~\ldots
$$
such that for some $k\in \N$, $g,h \in (0,1)$
$$
p_{2l}~=~\left\{ \begin{array}{l}
g ~~~~~\hbox{if~ } l=kj \hbox{~~~and }j \hbox{ is odd}  \\
h ~~~~~\hbox{if~ } l=kj \hbox{~~~and }j \hbox{ is even}  \\
{1\over 2} ~~~~~\hbox{else}~,
\end{array}\right.
$$
then 
\begin{equation}
\Phi^* (z) ~=~ {1\over 2h} {(1-2h)T_k^2(z) +(h-g)-\sqrt{
(T_k^2(z)-\k)^2-4\m}\over U_{k-1}(z)T_k(z)(1-z^2)}
\label{41b}
\end{equation}
where
\begin{equation}
\left\{ \begin{array}{l}
\k= g(1-h)+h(1-g) \\
\\
\m =g(1-g)h(1-h)
\end{array}\right.
\label{42}
\end{equation}
and the sign of the square root is defined such
that
\begin{equation}
\left|{T_k^2(z)-\k \over 2\sqrt{\m}}~+~\sqrt{{(T_k^2(z)-\k)^2
\over 4\m}-1}\right|~>~1~.
\label{42a}
\end{equation}
}

\bigskip

{\bf Proof.}~~Similary as in the first part
of the proof of Theorem 3.1 a contraction such that the 
convergent of the resulting continued fraction attains
successively the values of the $(k-1)$th, $(2k-1)$th,
$(3k-1)$th, $\ldots $ convergent, yields
\begin{equation}
\Phi ^*(z)~=~ \left[z-{1\over 2}\Psi (z)\right]^{-1}
\label{b3a}
\end{equation}
where 
\begin{eqnarray*}
\Psi (z)~&=&~
{~1~|\over |~z~}-\underbrace{{~{1\over 4}|\over|~z~}-
\ldots  -{~{1\over 4}|\over|~z~}}_{k-2} - {{1\over 2}g|\over |~z~}
- {{1\over 2}(1-g)|\over |~~~~z~~~~}\\
&&~~~~~~~~~~~~~~~-~
\underbrace{{~{1\over 4}|\over|~z~}-\ldots  -{~{1\over 4}|\over|~z~}}_{k-2}
- {{1\over 2}h|\over |~z~}-{{1\over 2}(1-h)|\over |~~~~z~~~~}
~-  {~{1\over 4}| \over |~z~} -~\ldots \\
&&\\
&=&~{2\over U_{k-1}(z)} \left[~U_{k-2}(z)
+ {~~g~~|\over |T_k(z)}-{h(1-g)|\over|~~T_k(z)~}-{g(1-h)|\over|~
~T_k(z)~} -{h(1-g)|\over|~~T_k(z)~}~-~ \ldots ~\right] ~.
\end{eqnarray*}
By (\ref{b3a}) and a further even contraction it now follows
that
$$
\Psi (z) ~=~
U_{k-1}(z)T_k(z) \left[{~~~~~1~~~~~|\over |T_k^2(z)-g}
- {g(1-g)h|\over |T_k^2(z)-\k}-{~~~~~\m~~~~~|\over|T_k^2(z)-\k}
-{~~~~~\m~~~~~|\over |T_k^2(z)-\k}~ -~ \ldots ~\right] ~
$$
where $\m$ and $\k$ are defined in (\ref{42}). The assertion is now
a consequence from standard results noting that the distribution with density
${2\over \pi}\sqrt{1-x^2}$ ($x \in [-1,1]$) has the Stieltjes transform
$$
{2\over \pi}\int_{-1}^1{\sqrt{1-x^2}\over w-x}dx 
~=~{{1}|\over |w}
~-~{{1\over 4}|\over |w}~-
~{{1\over 4}|\over |w}~-~\ldots ~=~2(w-\sqrt{w^2-1})
$$
where the square root is such that $|w+\sqrt{w^2-1}|>1$
(see e.g. VanAssche [14, p. 176]) \hfill \bull

\bigskip

{\bf Remark 4.2.}~~Note that the case $g=h$ was already
discussed in [10] where sieved
orthogonal polynomials on several intervals are
 considered (using the notation $g=h={1\over c+1}$). A delicate
analysis depending on (\ref{42a}) shows that 
the function in (\ref{41b})
has no poles at the zeros of the polynomial 
$(1-z^2)U_{k-1}(z)$ provided that  $g+h\leq 1$.
 If $g+h>1$ there are simple poles 
at $\pm 1$ with residues ${g+h-1\over 2hk}$ and  simple poles
at the zeros of $U_{k-1}(x)$ with residues ${g+h-1\over
hk}$. Moreover, it can be shown that $\Phi^* (z)$ has
no poles at the zeros of $T_k(z)$ if
$h\leq g$ and simple poles with residues ${h-g\over hk}$
if $h>g$. Finally we mention that  by the
 Perron-Stieltjes inversion formula the
absolute continuous component of $d\rho (x)$ is supported
on the set
$$
\overline{E_k (\k, \m)} ~:=~\overline{\left\{~
x\in [-1,1]~|~|T_k^2(x)-\k|^2<4\m \right\}}
$$
with density 
\begin{equation}
{d\rho (x)\over dx}~=~{1 \over 2\pi h}{\sqrt{
4\m-(T_k^2(x)-\k)^2}\over |U_{k-1}(x)||T_k(x)|(1-x^2)}
~~~~~x\in \overline{E_k (\k,\m)}
\label{44}
\end{equation}
where $\k$ and $\m$ are defined in (\ref{42}).

\bigskip

{\bf Theorem 4.3.}~~{\it Let ~$\lim_{l\to \infty}{\a_l\over l}
=a\geq 0$~ and ~$\lim_{l\to \infty}{\g_l\over l}=c \geq 0~$,
then $N_l^{(\as_l,\gs_l)}$ defined by (\ref{41a})
and (\ref{astern}) satisfies
\begin{equation}~~~~~~~
\lim_{l\to \infty}
{1\over l}N_l^{(\as_l,\gs_l)}(y)={2a+c+2 \over 2\pi }
\int_{-1}^y {\sqrt{
4\mac-(T_k^2(x)-\kac)^2}\over |U_{k-1}(x)||T_k(x)|(1-x^2)}
I\{ x\in \overline{E_k (\kac,\mac)}\}dx
\label{4b0}
\end{equation}
with
\begin{eqnarray}
\mac~&=&~\displaystyle{(2a+1)(c+1)(2a+c+1)\over (2a+c+2)^4} \nn \\
&&\\
\kac~&=&~\displaystyle{(2a+1)(c+2)+c(c+1) \over 
(2a+c+2)^2}\nn
\label{kappa}
\end{eqnarray}
}

\bigskip

{\bf Proof.}~~Because $\as _l = \a _{\lfloor (l-1)/k\rfloor}$,
 $\gs _l = \g _{\lfloor (l-1)/k\rfloor}$ are independent of $l\in
 \{nk+1,\ldots ,(n+1)k\}$, the zeros of
 $$
 ~C_{nk+1}^{(\as_{nk+1},\gs_{nk+1})}(x,k)~, 
 ~C_{nk+2}^{(\as_{nk+2},\gs_{nk+2})}(x,k)~, ~\ldots ,~ 
 ~C_{(n+1)k}^{(\as_{(n+1)k},\gs_{(n+1)k})}(x,k)~
 $$
are interlacing and consequently it is sufficient to prove
(\ref{4b0}) for the subsequence $l=nk$. 
The polynomial $C_{nk}^{(\as_{nk},\gs_{nk})}(x,k)$
is proportional to the polynomial in the denominator of the Stieltjes
transform of the probability measure $\xi_{\as_{nk},\gs_{nk},nk-1}^R$
which is obtained by terminating the continued fraction expansion 
for the Stieltjes transform of $\xi_{\as_{nk},\gs_{nk}}$ 
with density (\ref{w1}) and reversing the
corresponding sequence of $p_{2i}$'s (see (\ref{term})
and (\ref{termrev})). By the results
of Section 3 ~$\xi_{\as_{nk},\gs_{nk},nk-1}^R$ is "nearly" a uniform
distribution on the set  $\{x|~
C_{nk}^{(\as_{nk},\gs_{nk})}(x,k)=0\}$
and we have to distinguish the following two cases:

\medskip

{\bf A:)} $n=2m$, $l=2mk$: In this case it follows
from Remark 3.3 that $\xi_{\as_{2mk},\gs_{2mk},2mk-1}^R$ has equal masses at the zeros
of  $C_{2mk}^{(\as_{2mk},\gs_{2mk})}(x,k)$ and the
$p^{(2m)}_{2i}$ in the corresponding continued fraction expansion
$$
~~~~~~\int_{-1}^1 {d\xi_{\as_{2mk},\gs_{2mk},2mk-1}^R (x)\over z - x} ~=~
{~1~|\over |~z~}-{p_2^{(2m)}|\over|~~z~~}-{q_2^{(2m)}p_4^{(2m)}|\over|
~~~~~z~~~~~} 
- \ldots -{q_{4mk-4}^{(2m)}p_{4mk-2}^{(2m)}|\over|~~~~~~~z~~~~~~~} \ldots
$$
are given by $p_{2j}^{(2m)} =p_{2(2mk-j)}$ ($j=1,\ldots ,2mk-1$). 
From $p_{2kj}=A_{j-1}^{(2m-1)}$ ($j=1,\ldots ,2m-1$),
 (\ref{41}), and (\ref{lab3}) we have 
$$
p_{2kj}^{(2m)}~=~\left\{\begin{array}{l}
\displaystyle{\g_{2m-1}+2m-j\over 2\a_{2m-1}+\g_{2m-1}+
2(2m-j)}~~~\hbox{if }
j\hbox{ is odd}\\
\\
\displaystyle{2m-j\over 2\a_{2m-1}+\g_{2m-1}+2(2m-j)}~~~\hbox{if }
j\hbox{ is even}
\end{array}\right.
$$
($j\leq 2m-1$) and $p_{2l}^{(2m)} ={1\over 2}$ whenever
$l\neq jk$ ($l\leq 2mk-1$). From the assumption
of Theorem 4.3 it follows that 
\begin{equation} 
\lim_{m\to \infty}
p_{2l}^{(2m)}~=~p_{2l}^*~=~
\left\{\begin{array}{l}
\displaystyle{1\over 2} ~~~~~~~~~~~~~~~~~~~~~~~
~~~~~~\hbox{if }  l\neq jk \\
\\
\displaystyle{c+1\over 2a+c+2}~~~\hbox{if }
l=jk,~j\hbox{ is odd}\\
\\
\displaystyle{1\over 2a+c+2}~~~\hbox{if }
l=jk,~j\hbox{ is even}
\end{array}\right.
\label{limit}
\end{equation}
and by the same reasoning as in [7] we obtain
that the uniform distribution on the set 
$\{x|~C_{2mk}^{(\as_{2mk},\gs_{2mk})}(x,k) =0\}$
converges weakly to the distribution $\xi^*$
corresponding to the sequence $(p_2^*,p_4^*,p_6^*,\ldots )$
 in (\ref{limit}).
An application of Lemma 4.1 and Remark 4.2 shows that
$\xi ^*$ is absolute continuous with density 
given in (\ref{44}) where
$$
g={c+1\over 2a+c+2}~,~~h={1\over 2a+c+2}~,~~g+h\leq 1~,~~h\leq g~.
$$
This proves the assertion for the subsequence $l=2mk$.

\medskip

{\bf B:)}~~$n=2m+1$, $l=(2m+1)k$: By a similar argument as in
part A) we obtain that the probability measure $\xi_{\as_{(2m+1)k},\gs_{(2m+1)k},
(2m+1)k-1}^R$ converges
weakly to the distribution $\xi ^{**}$ corresponding to the sequence
$(p_2^{**},p_4^{**},p_6^{**},\ldots )$ where
$$
p_{2l}^{**}~=~
\left\{\begin{array}{l}
\displaystyle{1\over 2}  ~~~~~~~~~~~~~~~~~~~~~~~
~~~~~~\hbox{if }  l\neq jk \\
\\
\displaystyle{1\over 2a+c+2}~~~\hbox{if }
l=jk,~j\hbox{ is odd}\\
\\
\displaystyle{c+1\over 2a+c+2}~~~\hbox{if }
l=jk,~j\hbox{ is even}
\end{array}\right.~~.
$$
By Remark 4.2  (with $g={1\over 2a+c+2}$, $
h={c+1\over 2a+c+2}$, $g+h \leq 1$, $h\geq g$)
the absolute continuos component of $\xi ^{**}$
has the density
$$
f_{a,c}(x) ~=~
{2a+c+2 \over 2\pi (c+1)}{\sqrt{
4\mac-(T_k^2(x)-\kac)^2}\over |U_{k-1}(x)||T_k(x)|(1-x^2)}
~~~~~x\in \overline{E_k (\kac,\mac)}
$$
where $\kac, \mac$ are defined in (\ref{kappa})
(note that $\k$ and $\m$ defined by (\ref{42})
are symmetric in $g$ and $h$).
Moreover, there also exists a discrete part of $\xi^{**}$ 
at the zeros of $T_k(z)$ with equal masses
$$
{h-g\over kh}~=~{c\over k(c+1)}
$$
(see Remark 4.2). Now Theorem 3.1 (and its proof) show that
$\xi_{\as_{(2m+1)k},\gs_{(2m+1)k},(2m+1)k-1}^R$ has equal masses
$$
{\g_{2m}+1\over (2m+\g_{2m}+1)k}
$$
at the zeros of $T_k(x)$ and masses
$$
{1\over (2m+\g_{2m}+1)k} 
$$
at the remaining $2mk$ support points. This yields
\begin{eqnarray*}
&&{1\over (2m+1)k}N_{(2m+1)k}^{(\as_{(2m+1)k},\gs_{(2m+1)k})}(y)\\
&&\\
&&~=~ {(2m+\g_{2m}+1)k \over (2m+1)k}
\left[ {\g _{2m} + 1 \over (2m+\g_{2m}+1)k}
\# \{x\leq y|~T_k(x)=0\} \right.\\
&&  \hbox{\hskip 2.5cm}~+~{1\over (2m+\g_{2m}+1)k}\# \{ x\leq y|
~C_{(2m+1)k}^{(\as_{(2m+1)k},\gs_{(2m+1)k})}(x,k)=0,~T_k(x) \neq 0 \}\\
&&\left. \hbox{\hskip 2.5cm}~-~{\g _{2m}\over (2m+\g_{2m}+1)k}
\# \{ x\leq y|~T_k(x)=0\} \right] \\
&& \\
&&~=~{(2m+\g_{2m}+1) \over (2m+1)}\Biggl[
\int_{-1}^y d\xi_{\as_{(2m+1)k},\gs_{(2m+1)k},(2m+1)k-1}^R (x) 
\Biggr. \\
&& \Biggl. \hbox{\hskip 2.5cm}~-~{\g_{2m}\over (2m+\g_{2m}+1)k}
\# \{ x\leq y|~T_k(x)=0\} \Biggr] \\
&& \\
&& \longrightarrow
(c+1) \left[ \int_{-1}^y  d\xi ^{**}(x)
~-~ {c\over (c+1)k}
\# \{ x\leq y|~T_k(x)=0\} \right] \\
&& \\
&&=~{2a+c+2 \over 2\pi }
\int_{-1}^y {\sqrt{
4\m_{a,c} -(T_k^2(x)-\k_{a,c})^2}\over |U_{k-1}(x)||T_k(x)|(1-x^2)}
I\{ x\in \overline{E_k(\kac,\mac)}\}dx
\end{eqnarray*}
This proves (\ref{4b0}) for the subsequence
$l=(2m+1)k$ and by part A) for the sequence $l=nk$.
The assertion of the Theorem follows from the
 discussion at the beginning of this proof.
\hfill \bull

\bigskip

The following result states the analogue for the
asymptotic distribution of the zeros of
the \gsu polynomials of the second kind and
can be proved by similar arguments as given in the 
proof of Theorem 4.3.

\bigskip

{\bf Theorem 4.4.}~~{\it For $y\in [-1,1]$ let
$$
M_l^{(\as_l,\gs_l)}(y)~:=~\#\{x\leq y~|~\bags (x,k)=0~\}
$$
denote the number of zeros of 
$\bags (x,k)$ less or equal than $y$. If ~$\lim_{l\to \infty}
{\a_l\over l}=a\geq 0$~ and ~$\lim_{l\to \infty}
{\g_l\over l}=c \geq 0~$, then
$$
\lim_{l\to \infty}
{1\over l}M_l^{(\as_l,\gs_l)}(y)~=~{2a+c+2 \over 2\pi }
\int_{-1}^y {\sqrt{
4\m_{a,c}-(T_k^2(x)-\k_{a,c})^2}\over |U_{k-1}(x)||T_k(x)|(1-x^2)}
I\{ x\in \overline{E_k (\kac,\mac)}\}dx \nn
$$
where $\mac$, $\kac$ are defined in (\ref{kappa})}.

\bigskip

Note that in the case $c=0$ the limit distribution in Theorem
4.3 and 4.4 is exactly the distribution for which the
sieved Chebyshev polynomials of the first kind (introduced
in [10, p. 99] are orthogonal
with respect to. In the case considered here the discrete spectrum of this
measure is empty (the parameter $c$ in the lastnamed reference
corresponds to $2a+1$ which is obviously $\geq 1$ (see 
[10, p. 99]). The case $a=c=0$ and the trigonometric
identity yield the arcsin measure as limit distribution. \\
It is also worthwhile to mention  that an alternative proof  
of the asymptotic results for zeros of random walk polynomials 
$r_{nk}(x)=R_n(T_k(x)) - R_{n-2}(T_k(x))$ [or of the
second kind $s_{nk+k-1}(x)=U_{k-1}(x)R_n(T_k(x))$]
could be obtained if the asymptotic 
distribution of the zeros of the polynomials $R_n(x)
-R_{n-2}(x)$ [or $R_n(x)$] is known. This proof requires
additionally a careful inspection of the
inverse branches of $T_k(x)$ (see also [9]).

\medskip

We will conclude this paper by discussing the corresponding
limit statement for the zeros of the generalized Hermite
polynomials $\hg _n(x)$ defined in Section 2.
As pointed out in the introduction this result
has already been proved in  [8] as an
application of strong asymptotics for the generalized
Laguerre polynomials and the relation between 
these and the generalized Hermite polynomials
(see  [6, p. 156]). The proof given here 
is based on the characterizing property (\ref{ca1}) 
of the generalized Hermite polynomials (note that 
parametrization and standardization in [8] is different).

\bigskip

{\bf Theorem 4.5.}~~{\it For $y \in \R$ let
$$
N_l^{(\g_l)}(y)~=~\#\{~x\leq y~|~H_l^{(\g_l)}(x)=0~\}
$$ 
denote the number of zeros of $H_l^{(\g_l)}(x)$ less or
equal than $y$. If ~$\lim_{l \to \infty} {\g_l\over l}=c\geq 0$,
then
\begin{equation}
\lim_{l \to \infty} {1\over l}N_l^{(\g_l)}(\sqrt{l}y)~=~
{1\over \pi}\int_{-\infty}^\xi
{\sqrt{(c+1)-|x^2-(c/2)-1)|^2} \over |x|}
I\{x \in \overline{E_1}\}dx ~
\label{dens}
\end{equation}
where ${\overline E_1}=\{ x \in \R ~|~|x^2-(c/2)-1| 
<\sqrt{c+1}\}$.}

\bigskip

{\bf Proof.}~~We only consider the case of an
even subsequence $l=2m$, the odd case is treated
exactly in the same way as part B) in the proof of 
Theorem 4.3. \\
Let $\xi_{\g_l}$ denote the probability measure 
with density proportional to $|x|^{\g_l} \exp(-x^2)$
and let $\m_{2m-1}$ denote the uniform
distribution on the set
$$
\left\{ ~{x\over \sqrt{2m}}~|~H_{2m}^{(\g_{2m})}(x)=0~
\right\}~.
$$
From (\ref{ca1}), (\ref{rev}) and (\ref{canmom}) we obtain for the Stieltjes
transform of $\m_{2m-1}$
$$
\int_{-\infty}^\infty {d\m_{2m-1} (x)\over z - x} ~=~
~\int_{-\infty}^\infty {d\xi_{\g,2m-1}^R (x)\over z - 
 x/ \sqrt{2m}} ~=~
{~1~|\over |~z~}~-~{\at_1^{(2m)}|\over|~~z~~~}~-~ \ldots  ~-~
{\at_{2m-1}^{(2m)}|\over|~~~z~~~} 
$$
where
$$
\at_j^{(2m)}~=~\left\{ \begin{array}{l}
\displaystyle{{1\over 2}{2m-j\over 2m}}~~~~~~~~~~~~~~~
\hbox{if }j\hbox{ is even} \\
\\
\displaystyle{{1\over 2}{2m-j+\g_{2m}\over 2m}} ~~~~~~~\hbox{if }j
\hbox{ is odd} \\
\end{array} \right.~~.
$$
Thus we have 
\begin{equation}
\lim_{m \to \infty} \at_j^{(2m)}~=~\left\{ \begin{array}{l}
\displaystyle{{1\over 2}=: h}~~~~~~~~~~~~~~~
\hbox{if }j\hbox{ is even} \\
\\
\displaystyle{{1\over 2}(c+1)}=:g ~~~~~~\hbox{if }j\hbox{ is odd} \\
\end{array} \right.~
\label{last}
\end{equation}
and the 
the same arguments as in [7] show
that $\m_{2m-1}$ converges weakly to the
distribution $\m^*$ with Stieltjes transform
\begin{eqnarray*}
\Phi^*(z)~&=&~
{~1~|\over |~z~}~-~{~g~|\over|~z~}~ -~{~h~|\over|~z~}~ -~{~g~|\over|~z~}~ -~
{~h~|\over|~z~}~-~ \ldots  \\
&&\\
&=&~{~~~z~~|\over |z^2-g}-{~~~~~gh~~~~|\over |z^2 -(g+h)}
-{~~~~~gh~~~~|\over |z^2 -(g+h)}
-{~~~~~gh~~~~|\over |z^2 -(g+h)}~-\ldots \\
&&\\
&=&~
{z^2+h-g-\sqrt{[z^2-(g+h)]^2-4hg}\over 2hz}~.
\end{eqnarray*}
Here the square root is defined such that 
$$
\left| {z^2-(g+h)\over 2\sqrt{gh}}+\sqrt{{[z^2-(g+h)]^2
\over 4gh}-1}\right|~>~1
$$
and the
 second equality follows from the first by an even
contraction while the second equality implies
the third by the
same arguments as given in the proof of Lemma 4.1.
An application of (\ref{last}) (which implies $h\leq g$)
and the Perron-Stieltjes inversion formula 
show that $\m^*$ is absolute continuous with the density in
(\ref{dens}). This proves the assertion of the theorem.

\bigskip

{\bf Acknowledgements.}~~Parts of this paper were written while
the author was visiting the University of G\"ottingen. The author 
would like to thank the Institut f\"ur Mathematische Stochastik
for its hospitality and Dick Askey and Walter VanAssche for their help
with the references.

\bigskip

\begin{center}
REFERENCES\\
\end{center}
1. Abamowitz, M. and Stegun, I. (1964). {\it Handbook of Mathematical 
Functions}. Dover, New York \\
2. Al-Salam, W.,  Allaway, W.R. and Askey, R. (1984). Sieved
ultraspherical polynomials, {\it Trans. Amer. Math. Soc.},
{\bf 284}, 39--55. \\
3. Al-Salam, W. (1990). Characterization theorems for
orthogonal polynomials. In: {\it Orthogonal polynomials: Theory 
and Practice} (P. Nevai, Ed.), 1--24, NATO ASI Series
C 294, Kluwer, Dordrecht.\\
4. Charris, J.A. and Ismail, M.E.H. (1986). On sieved orthogonal
polynomials II: Random walk polynomials. {\it Canad. J. Math.}
{\bf 38}, 397--415.\\
5. Charris, J.A. and Ismail, M.E.H. (1993). Sieved orthogonal
polynomials VII: Generalized polynomial mappings. {\it Trans. Amer. Math.
Soc.} {\bf 340}, 71--93.\\
6. Chihara, T.S. (1978).
{\it Introduction to Orthogonal Polynomials}, Gordon \& Breach, New
York. \\
7. Dette, H. and Studden, W.J. (1992).
On a new characterization of the classical
orthogonal polynomials, {\it J. Approx. Theory}.,
{\bf 71}, 3--17.\\
8.  Gawronski, W. (1993). Strong asymptotics and asymptotic zero
distribution of
Laguerre polynomials $L_n^{(an+\a)}(x)$ and Hermite polynomials
$H_n^{(an+\a)}(x)$. Analysis, {\bf 13}, 29--68.\\
9. Geronimo, J.S. and VanAssche, W. (1988). Orthogonal polynomials on several
intervals via a polynomial mapping. {\it Trans.  Amer. Math. Soc.}
{\bf 308}, 559--581. \\
10. Ismail, M.E.H. (1986). On sieved orthogonal polynomials III:
Orthogonality on several intervals. {\it Trans. Amer. Math. Soc.}
{\bf 294}, 89--111. \\
11. Perron, O. (1954).
{\it Die Lehre von den Kettenbr\"uchen (Band I and II)}, B.G. Teubner,
Stuttgart. \\
12. Rogers, L.J. (1895). Third memoir on the expansion
of certain infinite products. {\it Proc. London Math. Soc.} {\bf 26},
15--32. \\
13. Szeg\"o, G. (1975). {\it Orthogonal polynomials}. American Mathematical
Society Colloquium Publications, Vol. {\bf 23}, Providence, RI. \\
14. VanAssche, W. (1987). {\it Asymptotics for orthogonal polynomials}.
Lecture Notes in Mathematics, No. 1265, Springer Verlag, Berlin, New York.\\
15. Wall, H.S. (1948).
{\it Analytic theory of continued fractions}, Van Nostrand, New York.

\end{document}